# A Faster Product for $\pi$ and a New Integral for $\ln\dfrac{\pi}{2}$

## Jonathan Sondow

**1. INTRODUCTION.** In [**15**] we derived an infinite product representation of $e^\gamma$, where $\gamma$ is Euler's constant:

$$e^\gamma = \left(\frac{2}{1}\right)^{1/2}\left(\frac{2^2}{1\cdot 3}\right)^{1/3}\left(\frac{2^3\cdot 4}{1\cdot 3^3}\right)^{1/4}\left(\frac{2^4\cdot 4^4}{1\cdot 3^6\cdot 5}\right)^{1/5}\cdots. \qquad (1)$$

Here the $n$th factor is the $(n+1)$th root of the product

$$\prod_{k=0}^{n}(k+1)^{(-1)^{k+1}\binom{n}{k}}.$$

In the process we noticed a strikingly similar product representation of $\pi$:

$$\frac{\pi}{2} = \left(\frac{2}{1}\right)^{1/2}\left(\frac{2^2}{1\cdot 3}\right)^{1/4}\left(\frac{2^3\cdot 4}{1\cdot 3^3}\right)^{1/8}\left(\frac{2^4\cdot 4^4}{1\cdot 3^6\cdot 5}\right)^{1/16}\cdots. \qquad (2)$$

In this note we give three proofs of (2). The third leads to an analog for $\ln(\pi/2)$ of integrals for $\ln(4/\pi)$ [**14**] and $\gamma$ [**13**], [**14**], [**15**]:

$$\ln\frac{\pi}{2} = -\int_{[0,1]} \frac{1-x}{(1+x)\ln x}\,dx, \qquad (3)$$

$$\ln\frac{4}{\pi} = -\iint_{[0,1]^2} \frac{1-x}{(1+xy)\ln xy}\,dx\,dy, \qquad (4)$$

$$\gamma = -\iint_{[0,1]^2} \frac{1-x}{(1-xy)\ln xy}\,dx\,dy.$$

Using (3), we sketch a derivation of (1) and (2) from the same function (a form of the polylogarithm [**7**]), accounting for the resemblance between the two products. The function also leads to a product for $e$ (due to J. Guillera [**5**]),

$$e = \left(\frac{2}{1}\right)^{1/1}\left(\frac{2^2}{1\cdot 3}\right)^{1/2}\left(\frac{2^3\cdot 4}{1\cdot 3^3}\right)^{1/3}\left(\frac{2^4\cdot 4^4}{1\cdot 3^6\cdot 5}\right)^{1/4}\cdots, \qquad (5)$$

surprisingly close to product (1) for $e^{\gamma}$.

**2. THE ALTERNATING ZETA FUNCTION.** The logarithm of product (1), namely,

$$\gamma = \sum_{n=1}^{\infty} \frac{1}{n+1} \sum_{k=0}^{n} (-1)^{k+1} \binom{n}{k} \ln(k+1), \qquad (6)$$

reminded us of the series (see [**6**] and [**11**])

$$\zeta^*(s) = \sum_{n=0}^{\infty} \frac{1}{2^{n+1}} \sum_{k=0}^{n} (-1)^k \binom{n}{k} (k+1)^{-s} \quad (s \in \mathbf{C}), \qquad (7)$$

which gives the analytic continuation of the alternating zeta function $\zeta^*(s)$. The latter is defined by the Dirichlet series (see [**12**])

$$\zeta^*(s) = \sum_{k=1}^{\infty} \frac{(-1)^{k-1}}{k^s} \quad (\Re(s) > 0). \qquad (8)$$

(For example, using the classic formula $\zeta^*(1) = \ln 2$ for the alternating harmonic series—for a new proof see [**12**]—one can derive the series $\ln 2 = \sum_{n \geq 1} (2^n n)^{-1}$ from (7) by considering it when $s = 1$.) Differentiating (7) termwise and substituting the value of the derivative of $\zeta^*$ at $s = 0$,

$$\zeta^{*\prime}(0) = \frac{1}{2} \ln \frac{\pi}{2} \qquad (9)$$

(see [**11**]), yields the series

$$\ln \frac{\pi}{2} = \sum_{n=1}^{\infty} \frac{1}{2^n} \sum_{k=0}^{n} (-1)^{k+1} \binom{n}{k} \ln(k+1), \qquad (10)$$

and exponentiation produces product (2).

**3. WALLIS'S PRODUCT AND EULER'S TRANSFORM.** The pair of infinite products (1) and (2) calls to mind another pair, Wallis's product for $\pi$ [**17**] and Pippenger's product for $e$ [**10**]:

$$\frac{\pi}{2} = \frac{2}{1}\frac{2}{3}\frac{4}{3}\frac{4}{5}\frac{6}{5}\frac{6}{7}\frac{8}{7}\cdots, \tag{11}$$

$$\frac{e}{2} = \left(\frac{2}{1}\right)^{1/2}\left(\frac{2}{3}\frac{4}{3}\right)^{1/4}\left(\frac{4}{5}\frac{6}{5}\frac{6}{7}\frac{8}{7}\right)^{1/8}\cdots. \tag{12}$$

It is interesting to note that products (2) and (12), whose factors have exponents $1/2^n$, converge rapidly to numbers $\pi/2$ and $e/2$ whose irrationality has been proved (see, for example, [**9**]), whereas product (1), with exponents $1/(n+1)$, converges less rapidly to a number $e^\gamma$ whose (expected) irrationality has not yet been proved.

We give a second proof of (2), using (11) and Euler's transformation of series

$$\sum_{n=1}^{\infty}(-1)^{n-1}a_n = \sum_{n=0}^{\infty}\frac{1}{2^{n+1}}\sum_{k=0}^{n}(-1)^k\binom{n}{k}a_{k+1}, \tag{13}$$

valid for any convergent series of complex numbers [**7**, sec. 33B], [**11**]. Applying (13) to the logarithm of Wallis's product

$$\ln\frac{\pi}{2} = \sum_{n=1}^{\infty}(-1)^{n-1}\ln\frac{n+1}{n} \tag{14}$$

gives

$$\ln\frac{\pi}{2} = \sum_{n=0}^{\infty}\frac{1}{2^{n+1}}\sum_{k=0}^{n}(-1)^k\binom{n}{k}\ln\frac{k+2}{k+1}. \tag{15}$$

If we replace $n$ by $n-1$, write the last logarithm as $\ln(k+2) - \ln(k+1)$, and the sum on $k$ as the difference of two sums in the first of which we replace $k$ by $k-1$, then the recursion $\binom{n-1}{k-1} + \binom{n-1}{k} = \binom{n}{k}$ leads to (10), completing the second proof of (2). The first proof is basically the same, because in [**11**] we use Wallis's product to evaluate (9), and we take the Euler transform of (8) to get (7) for complex $s$ with $\Re(s) > 0$.

Products (12) and (11) are linked by Stirling's asymptotic formula $n! \sim (n/e)^n\sqrt{2\pi n}$: the formula is proved in [**2**] using (11) and is used in [**10**] to establish (12). Products (1) and (2) are linked by transformations: a hypergeometric one [**15**] for (1) and Euler's for (2). (To strengthen the link, we can write series (14) and (15) as integrals of hypergeometric functions—compare [**15**, Proof 1]—and then obtain (15) from (14) by a hypergeometric transformation equivalent to (13).) However, this link does not explain the remarkable resemblance between (1) and (2).

Euler's transformation accelerates the rate of convergence of a slowly converging series such as (14) (see [**7**, sec. 35B]). Thus, product (2) converges faster than product (11), as Figure 1 shows.

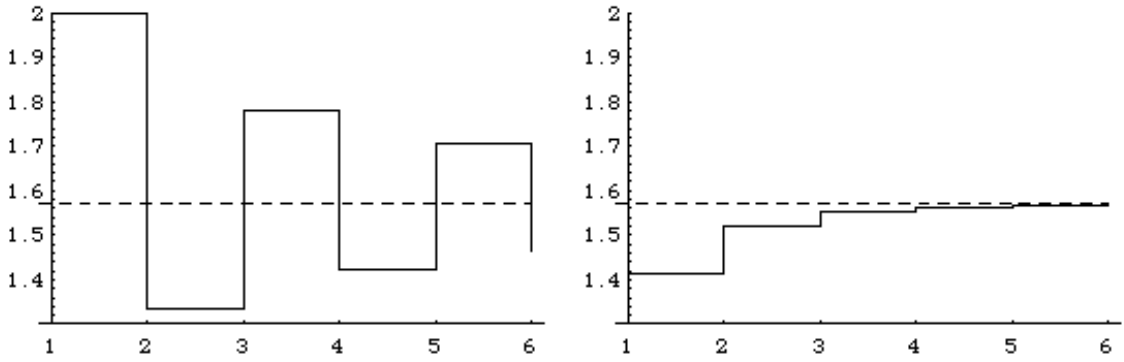

**Figure 1.** Partial products of Wallis's product and its Euler transform

**4. AVOIDING EULER.** A third proof of (2) (due in part to S. Zlobin [**18**]) avoids using Euler's transformation altogether (compare the proofs avoiding hypergeometric functions in [**15**]). We show that

$$I := \int_0^\infty \int_0^1 \sum_{n=0}^\infty x^y \left(\frac{1-x}{2}\right)^{n+1} dx\, dy = \ln \frac{\pi}{2}. \tag{16}$$

This implies (2), because if we factor $(1-x)$ from the integrand and use the binomial theorem, then termwise integration (justified since the integrand is majorized by the series $\sum 2^{-n-1}$) yields (15) and, therefore, (2). To prove (16), we use the geometric series summations

$$\sum_{n=0}^\infty \left(\frac{1-x}{2}\right)^{n+1} = \frac{1-x}{1+x} = \frac{(1-x)^2}{1-x^2} = \sum_{n=0}^\infty (1-x)^2 x^{2n} \tag{17}$$

to write

$$I = \int_0^\infty \int_0^1 \sum_{n=0}^\infty (1-x)^2 x^{y+2n}\, dx\, dy.$$

The integrand is majorized by $\sum (n+1)^{-2}$ (because

$$\max_{0 \le x \le 1} (1-x)^2 x^{2n} = \left(\frac{1}{n+1}\right)^2 \left(\frac{n}{n+1}\right)^{2n} < \frac{1}{(n+1)^2}$$

and $x^y \le 1$), so we may perform the integrations term by term, which by invoking (11) gives

$$I = \sum_{n=0}^{\infty} \ln \frac{(2n+2)^2}{(2n+1)(2n+3)} = \ln \frac{2^2}{1 \cdot 3} \frac{4^2}{3 \cdot 5} \frac{6^2}{5 \cdot 7} \cdots = \ln \frac{\pi}{2}.$$

This proves (16) and completes the third proof of (2).

*Proof of* (3). Equation (16) and the first equality in (17) yield

$$\int_0^\infty \int_0^1 x^y \frac{1-x}{1+x} \, dx \, dy = \ln \frac{\pi}{2}.$$

Reversing the order of integration (permitted since the integrand is nonnegative), we integrate with respect to $y$ and arrive at formula (3).

Alternatively, one can derive (3) from (4) by making the change of variables $u = xy, v = 1 - x$ and integrating with respect to $v$: the result is $\ln 2$ minus integral (3) (with $u$ in place of $x$), and equality (3) follows.

**5. RELATING THE PRODUCTS FOR $\pi$ AND $e^\gamma$.** Recall that we derived product (2) for $\pi$ from the alternating zeta function $\zeta^*(s)$. Omitting details, we sketch a derivation of product (1) for $e^\gamma$ from a generalization of $\zeta^*(s)$. This accounts for the resemblance between the two products. (Formulas (18) and (19) are due to J. Guillera [**5**].)

We generalize series (7) for $\zeta^*(s)$ by defining the function

$$f(t,s) = \sum_{n=0}^{\infty} t^{n+1} \sum_{k=0}^{n} (-1)^k \binom{n}{k} (k+1)^{-s} \quad (-1 < t < 1, \ s \in \mathbf{C}), \tag{18}$$

so that $f(1/2, s) = \zeta^*(s)$. Using integral (16) but replacing $(1-x)/2$ with $t(1-x)$, we can show that the formula obtained from (3) and (9),

$$\zeta^{*\prime}(0) = -\frac{1}{2} \int_0^1 \frac{1-x}{(1+x) \ln x} \, dx,$$

extends to

$$f'(t,0) = -t^2 \int_0^1 \frac{1-x}{(1-t(1-x)) \ln x} \, dx, \tag{19}$$

where the prime $'$ is shorthand for $\partial/\partial s$.

We now derive product (1) by evaluating the integral $\int_0^1 t^{-1} f'(t,0) \, dt$ in two different ways. On the one hand, a glance at (18) reveals that this integral equals the right side of (6). On the other hand, substituting (19) into the integral and reversing the order of integration gives

$$\int_0^1 \frac{f'(t,0)}{t} \, dt = -\int_0^1 \int_0^1 \frac{t(1-x)}{(1-t(1-x)) \ln x} \, dt \, dx = \int_0^1 \left( \frac{1}{\ln x} + \frac{1}{1-x} \right) dx.$$

The last is a classical integral for Euler's constant [**1**, sec. 10.3], [**15**], and (6) follows, implying (1).

Other products can be derived in the same way. For example, exponentiating the integral $\int_0^1 t^{-2} f'(t,0)\,dt = 1$ gives product (5) for $e$, which converges more slowly than Pippenger's product for $e$, because of the exponents $1/n$ in (5), versus $1/2^n$ in (12).

In order to identify the function $f(t,s)$, we reverse the order of summation in (18) and sum the resulting series on $n$. We then replace $k$ with $k-1$, obtaining

$$f(t,s) = \sum_{k=0}^{\infty} \frac{(-1)^k}{(k+1)^s} \sum_{n=k}^{\infty} \binom{n}{k} t^{n+1} = -\sum_{k=1}^{\infty} \frac{\bigl(t/(t-1)\bigr)^k}{k^s} \tag{20}$$

for $t$ satisfying $-1 < t \leq 1/2$ and for suitable $s$. Therefore, $f(t,s)$ and $\zeta^*(s)$ are related to the function

$$F(t,s) = \sum_{k=1}^{\infty} \frac{t^k}{k^s} \quad (-1 \leq t < 1,\ \Re(s) > 0)$$

by the formulas

$$f(t,s) = -F(t/(t-1), s), \tag{21}$$

$$\zeta^*(s) = -F(-1, s),$$

for appropriate $t$ and $s$. (With $t = 1/2$, equations (18) and (20) verify that formulas (7) and (8) for $\zeta^*(s)$ agree.) The function $F(t,s)$, a special case of the Lerch zeta function $\Phi(z,s,v)$ (see [**3**, sec. 1.11], [**16**, Sec. 64]), is the polylogarithm $\text{Li}_s(t)$ when $s$ is an integer [**7**, p. 189], [**16**, secs. 25, 64]. Relations (18) and (21) lead to an analytic continuation of $F(t,s)$, and thus of the polylogarithm.

*209 West 97th Street, New York, NY 10025*
*jsondow@alumni.princeton.edu*